\documentclass{amsart}
\usepackage{amssymb,latexsym,amsmath,amscd,amsthm,amsfonts}

\begin{document}

\title [Unitary Representations] {Unitary representations of\\compact quantum groups}
\author{Raluca Dumitru}
\address{Department of Mathematical Sciences\\ University of Cincinnati\\ Cincinnati, OH}
\address{Institute of Mathematics of the Romanian Academy\\ Bucharest, Romania}
\maketitle

\ 

\begin{abstract}Let $v$ be the right regular representation of a compact quantum group $G$. Then (\cite{wor2}, Proposition 4.4) $v$ contains all irreducible representations of $G$ and each irreducible representation enters $v$ with the multiplicity equal to its dimension. The result is certainly known for classical compact groups (\cite{ross}). The aim of this paper is to give a short survey on this subject and to provide a different proof of Woronowicz's result. The proof is an adaptation of the corresponding result for classical compact groups and provides a concrete decomposition of the right regular representation in irreducible components.
\end{abstract}

\

\

\ Quantum groups have been studied in the last decades by a large number of mathematicians and physicists. There are several approaches to the theory of quantum groups but the one used in this paper is the C*-algebraic approach introduced by S. L. Woronowicz (\cite{wor1}, \cite{wor2}).
\newtheorem{definition}{Definition}
\begin{definition}[\cite{wor2},Definition 2.1] A compact quantum group $G$ is given by a pair ($A$,$\Delta$) with $A$ a unital C*-algebra and $\Delta : A \to A\otimes A$ a unital $\ast$-homomorphism with the following properties:
\begin{enumerate}
\item The diagram is commutative
\[\begin{CD} A @>\Delta>> A\otimes A\\@V\Delta VV @VV\Delta \otimes id V\\A\otimes A @>>id\otimes \Delta > A\otimes A\otimes A
\end{CD}
\]
\ 

\item \[ \overline{(A\otimes 1)\Delta (A)} = A\otimes A \]
\ 
\[ \overline{(1\otimes A)\Delta (A)} = A\otimes A \]

\end{enumerate}
\end{definition}
\

\ Note that the tensor product is the minimal tensor product and that there is no separability assumption on $A$ as in Woronowicz's original definition (\cite{wor2}, Definition 2.1). The separability assumption was used by Woronowicz for proving the existence of the Haar state. In \cite{van} Van Daele proved that this assumption is not necessary by showing that a Haar state exists and is unique without assuming that $A$ is separable. The result is the following:
\newtheorem{proposition}[definition]{Proposition}

\begin{proposition}[\cite{wor2}, Theorem 2.3; \cite{van}, Proposition 2.4] Let $G$=($A$,$\Delta$) be a compact quantum group. Then there exists a unique state $h$ on $A$ with the property
\[ (id\otimes h)\Delta (a)=(h\otimes id)\Delta (a)= h(a)\cdot 1_{A}\]  
\end{proposition}
\ Then $h$ is called the Haar state on $G$. Note that the Haar state is not always faithful as in the case of classical compact groups. There exists however a subset of $A$, denoted by $\mathcal{A}$, on which $h$ is faithful. $\mathcal{A}$ is called the Woronowicz-Hopf $\ast$-algebra and is defined in terms of irreducible unitary representations of $G$. Applying the GNS construction to $h$ we obtain a Hilbert space $H_{h}$, a representation $\pi_{h}$ of $A$ and a cyclic vector $\xi_{h}$.
\ If $a\in A$ (or $a\in \mathcal{A}$) and $f$ is a bounded linear functional (on $A$ or $\mathcal{A}$), we denote (\cite{wor1}, \cite{wor2})
\[
\ a \ast f = (f \otimes \iota)(\Delta(a))
\]
\[
\ f \ast a = (\iota \otimes f)(\Delta(a))
\]

\ Both $a \ast f$ and $f \ast a$ are elements of $A$ (respectively $\mathcal{A}$).

\ Let $H$ be a Hilbert space and denote by $K(H)$ the set of compact operators on $H$. A unitary representation $u$ of $G$ acting on a Hilbert space $H$ is a unitary element in the multiplier algebra $\mathcal{M} (K(H)\otimes A)$ with the property $(id\otimes \Delta)u=u_{12}u_{13}$. The leg numbering notation used ($u_{12}$, $u_{13}$) is the one obtained from the two possible inclusions of $K(H)\otimes A$ in $K(H)\otimes A \otimes A$. The representation $u$ is called irreducible if the operators intertwining $u$ (i.e. $S\in B(H_{u})$ such that $(S\otimes 1)u=u(S\otimes 1)$, where $H_{u}$ denotes the Hilbert space on which $u$ acts) are scalar multiples of the identity. If $u$ is irreducible then it is finite dimensional (i.e. $dimH_{u}<\infty$) and if we denote by $d=dimH_{u}$ then $d$ is called the dimension of $u$. By choosing an orthonormal basis in $H_{u}$, the Hilbert space of an irreducible representation $u$, we can identify $u$ with a matrix $(u_{ij})_{1\leqslant i,j\leqslant d}$ over $A$. 

\ Two unitary representations $u$ and $v$ are called equivalent (denoted by $u\sim v$) if and only if there exists $U:H_{u}\to H_{v}$ unitary such that $(U\otimes 1)u=v(U\otimes 1)$. Let $\widehat G$ denote the dual of $G$, i.e. the set of all unitary equivalence classes of irreducible representations of $G$ (\cite{wor1}, \cite{wor2}). For each $\alpha \in \widehat G$, denote by $\ u ^{\alpha}$ a representative of each class. Let $ \{u_{ij}^{\alpha}\}_{1\leqslant i,j\leqslant d_\alpha}\subset A$ be the matrix elements of $\ u ^{\alpha}$. Then the Woronowicz-Hopf $\ast$-algebra $\mathcal{A}$ is defined by
\[ \mathcal{A}=lin\{u^{\alpha}_{ij}\mid u^{\alpha}\in \widehat G,\: 1\leqslant i,j\leqslant d_{\alpha}\}\]

\ Then $\mathcal{A}$ is a dense $\ast$-subalgebra of $A$ (\cite{wor2}, Propositions 5.2, 6.1, 6.2).

\ For any $\alpha \in \widehat G$ one can construct the contragradient representation
\[ \overline{u^{\alpha}}=\sum_{i,j=1}^{d_{\alpha}}m_{ij}\otimes {u^{\alpha}_{ij}}^{\ast}\]

\ Note that $\overline{u^{\alpha}}$ is an irreducible representation of $G$ but is not necessarily unitary. However, it is equivalent to a unitary representation of $G$. Using (\cite{wor1}, Theorem 5.2) one can show that $\forall$ $\alpha \in \widehat G$, $\exists$ $\beta \in \widehat G$ such that $\overline{u^{\alpha}}$ is equivalent to $u^{\beta}$. The representation $u^{\beta}$ will be called the conjugate representation and will be denoted by $u^{\overline \alpha}$. Notice that due to this equivalence the coefficients of the conjugate representation are linear combinations of ${u_{ij}^{\alpha}}^{\ast}$ and both 
$lin\{u^{\overline{\alpha}}_{ij}\mid 1\leqslant i,j\leqslant d_{\overline{\alpha}}\}$ and $lin\{{u_{ij}^{\alpha}}^{\ast}\mid 1\leqslant i,j\leqslant d_{\alpha}\}$ are linear basis in $p_{\alpha}H_{h}$.  

\ Let $\alpha \in \widehat{G}$ and set $\chi_{\alpha}=\sum_{i=1}^{d_{\alpha}}u^{\alpha}_{ii}$. Then $\chi_{\alpha}$ is called the character of $u^{\alpha}$. If $H_{\alpha}$ is the finite dimensional Hilbert space of the representation $u^{\alpha}$, then by (\cite{wor1}, \cite{wor2}) there exists a unique positive invertible operator $F_{\alpha}\in B(H_{\alpha})$ that intertwines $u^{\alpha}$ with its double contragradient representation $\overline{\overline{u^{\alpha}}}$, such that $Tr(F_{\alpha})=Tr(F_{\alpha}^{-1})$. The operator $F_{\alpha}$ can be represented as a matrix $F_{\alpha}=[f_{1}(u^{\alpha}_{ij})]_{i,j}$, with $f_{1}$ linear functional on $\mathcal{A}$. Set $M_{\alpha}=Tr(F_{\alpha})$ and $a_{\alpha}=M_{\alpha}(f_{1}\ast \chi_{\alpha})^{\ast}$. Clearly, $a_{\alpha}\in \mathcal{A}$. 
\newtheorem{theorem}[definition]{Theorem}
\begin{theorem}[\cite{wor1}, Theorem 5.7]\label{comp} With the notations introduced above the following relations hold:
\begin{enumerate}
\item $h(u^{\alpha}_{mk}(u^{\beta}_{nl})^{\ast})=\frac{1}{M_{\alpha}}\delta_{\alpha \beta}\delta_{mn}f_{1}(u^{\alpha}_{lk})$
\item $h((u^{\alpha}_{km})^{\ast}u^{\beta}_{ln})=\frac{1}{M_{\alpha}}\delta_{\alpha \beta}\delta_{mn}f_{-1}(u^{\alpha}_{lk})$
\item $h(a_{\alpha}u^{\beta}_{nl})=\frac{1}{M_{\alpha}}\delta_{\alpha \beta}\delta_{nl}$
\item $h((u^{\beta}_{nl})^{\ast}(a_{\alpha})^{\ast})=\frac{1}{M_{\alpha}}\delta_{\alpha \beta}\delta_{nl},$
\end{enumerate}
for all $\alpha$, $\beta \in \widehat G$ and for all $k,l, m, n \in \{1,...,d_{\alpha}\}$, where $f_{-1}(u^{\alpha}_{ij})$ are the entries of the matrix $F_{\alpha}^{-1}$.
\end{theorem}

\ Let $u^{\alpha}\in K(H_{\alpha})\otimes A$ be an irreducible representation of $G$. Then $H_{\alpha}$ is finite dimensional and denote $d_{\alpha}=dimH_{\alpha}$. Define the following subsets of the Hopf algebra $\mathcal{A}$:
\[ \mathcal{A}_{u^{\alpha}}=lin\{u^{\alpha}_{ij}\mid \: 1\leqslant i,j\leqslant d_{\alpha}\}\]
\[ \mathcal{A}_{u^{\alpha}}^{k}=lin\{u^{\alpha}_{ki}\mid \: 1\leqslant i\leqslant d_{\alpha}\}\]
with $k\in\{1,...,d_{\alpha}\}$.
\ Let
\[H_{h}(u^{\alpha})=\overline{\pi_{h}(\mathcal{A}_{u^{\alpha}})\xi_{h}}^{\| \cdot\|}\subseteq H_{h}\] and
\[H_{h}^{k}(u^{\alpha})=\overline{\pi_{h}(\mathcal{A}_{u^{\alpha}}^{k})\xi_{h}}^{\| \cdot\|}\subseteq H_{h}(u^{\alpha})\subseteq H_{h}.\]

\ Since for $\alpha,\beta\in\widehat G$, $\alpha\ne\beta$, 
\[\langle\pi_{h}(u^{\alpha}_{ij})\xi_{h}, \pi_{h}(u^{\beta}_{kl})\xi_{h}\rangle=h({u^{\beta}_{kl}}^{\ast}u^{\alpha}_{ij})=0,\: \forall i,j\in\{1,...,d_{\alpha}\},\: \forall k,l\in\{1,...,d_{\beta}\}\]
then 
\begin{equation}\label{dirsum1}
H_{h}=\underset{\alpha\in \widehat G}{\sum}\oplus H_{h}(u^{\alpha}),
\end{equation} 
as a direct sum of orthogonal spaces. Also,
\begin{equation}\label{dirsum2} 
H_{h}(u^{\alpha})=\overset{d_{\alpha}}{\underset{k=1}{\sum}}\oplus H_{h}^{k}(u^{\alpha}),
\end{equation}
a direct sum of orthogonal spaces.

\ Among the unitary representations of $G$ a special role is played by the right regular representation, defined by S. L. Woronowicz as follows: 

\begin{theorem}[\cite{wor2}, Theorem 5.1]\label{reg} With the notation introduced above
\begin{enumerate}
\item There exists unique $v\in \mathcal{M}(K(H)\otimes A)$ such that
\[ [(id\otimes f)v]\pi_{h}(a)\xi_{h}=\pi_{h}(f \ast a)\xi_{h}\]
\ for any $f$ continuous linear functional on $A$, any $a\in A$.
\item $v$ is a unitary representation of $G$.
\item The set
\[ \{(\rho \otimes id)v\:\mid\:\rho \:continuous\: linear\: functional\: on\: K(H_{h})\}\]
is dense in $A$.

\end{enumerate}
\end{theorem}

\ Then $v$ is called the right regular representation of the compact quantum group $G$. The set of intertwiners of $v$ and the definition of $C^{\ast}(G)$ are given in the following theorem:

\begin{theorem}[\cite{wor2}, Theorem 4.1] \label{inter} Let $v$ be the right regular representation of a compact quantum group $G=(A,\Delta)$. Let $C^{\ast}(G)$ be the norm closure of the set of all operators of the form
\[ \mathcal{F}_{v}(a)=(id\otimes ha)(v^{\ast}),\]
where $a\in A$. Then 
\begin{enumerate}
\item $v\in \mathcal{M}(C^{\ast}(G)\otimes A)$
\item An operator $X\in B(H_{h})$ intertwines $v$ with itself (i.e. $(X\otimes id)v=v(X\otimes id)$) if and only if $X$ commutes with all elements of $C^{\ast}(G)$.
\end{enumerate}
\end{theorem}

\ Consider now $V\in B(H_{h}\otimes H_{h})$ defined by $V=(id\otimes \pi_{h})v$. The following result is implicitly contained in (\cite{wor1}, \cite{wor2}):

\newtheorem{lemma}[definition]{Lemma}
\begin{lemma}\label{V} $V(\pi_{h}(a)\xi_{h}\otimes \xi)=(\pi_{h}\otimes \pi_{h})(\Delta(a))(\xi_{h}\otimes \xi)$, $\forall a\in A$, $\forall \xi\in H_{h}$.
\end{lemma}

With $H_{h}(\overline{u^{\alpha}})$ introduced above, the orthogonal projection onto this subspace of $H_{h}$ is $p_{\alpha}=F_{v}(a_{\alpha})$. Using the multiplication in $C^{\ast}(G)$ as described in (\cite{wor2}, formula 4.4) one can show that $p_{\alpha}$ is a central projection in $C^{\ast}(G)$. By equation \ref{dirsum1} above, $\sum p_{\alpha}=Id_{B(H_{h})}$. Therefore 
\[C^{\ast}(G)=\overline{\sum_{\alpha \in \widehat{G}}C^{\ast}(G)p_{\alpha}},\]
norm closure in $B(H_{h})$, where $\sum_{\alpha \in \widehat{G}}C^{\ast}(G)p_{\alpha}$ denotes the set of all finite sums.

The structure of $C^{\ast}(G)p_{\alpha}$ is described by the next proposition:
\begin{proposition}\label{factor} For every $\alpha \in \widehat{G}$, there is an isomorphism between $C^{\ast}(G)p_{\alpha}\subseteq B(H_{h}(u^{\alpha}))$ and the full matrix algebra $\mathcal{M}_{d_{\alpha}\times d_{\alpha}}$.
\end{proposition}
\begin{proof} Let $\alpha \in \widehat{G}$. Set $c_{ij}=M_{\alpha}(u^{\alpha}_{ij}\ast f_{1})^{\ast}\in \mathcal{A}$, for all $1\leqslant i,j\leqslant d_{\alpha}$. Note that $c_{ij}$ are linearly independent and
\[ lin\{c_{ij} \mid 1\leqslant i,j\leqslant d_{\alpha}\}=\mathcal{A}_{\overline{u^{\alpha}}}.\]
For each $1\leqslant i,j\leqslant d_{\alpha}$ let $E_{ij}=F_{v}(c_{ij})^{\ast}F_{v}(a_{\alpha})$. Then $E_{ij}(u^{\alpha}_{rs})^{\ast}=\delta_{sj}(u^{\alpha}_{ri})^{\ast}$ and therefore $E_{kl}E_{ij}=\delta_{il}E_{kj}$. The operator $E_{ij}$ can be represented as a block diagonal matrix. Each $d_{\alpha}\times d_{\alpha}$ block on the diagonal is given by the matrix $m_{ij}=[a_{kl}]_{1\leqslant k,l\leqslant d_{\alpha}}$, for all $1\leqslant i,j\leqslant d_{\alpha}$, where 
\[a_{kl}=
\begin{cases}
1, &\text{if k=i and l=j}\\
0, &\text{otherwise.}
\end{cases}\]

 The application $E_{ij}\longmapsto m_{ij}$, ${1\leqslant i,j\leqslant d_{\alpha}}$ establishes an isomorphism between $C^{\ast}(G)p_{\alpha}$ and the full matrix algebra $\mathcal{M}_{d_{\alpha}\times d_{\alpha}}$. Since the algebras $C^{\ast}(G)p_{\alpha}$ and $\mathcal{M}_{d_{\alpha}\times d_{\alpha}}$ are finite dimensional, they are $\ast$-isomorphic.
\end{proof}

\newtheorem{remark}[definition]{Remark}
\begin{remark} Since the unit $p_{\alpha}$ of $B(H_{h}(u^{\alpha}))$ belongs to $C^{\ast}(G)p_{\alpha}$, then $C^{\ast}(G)p_{\alpha}$ is a subfactor of $B(H_{h}(u^{\alpha}))$. Using (\cite{jacob}, page 118) $B(H_{h}(u^{\alpha}))$ is isomorphic to $C^{\ast}(G)p_{\alpha}\otimes (C^{\ast}(G)p_{\alpha})'$ and, by reasons of dimension, $(C^{\ast}(G)p_{\alpha})'$ is isomorphic to the full matrix algebra $\mathcal{M}_{d_{\alpha}\times d_{\alpha}}$.
\end{remark}
\ The description of the right regular representation is given by S. L. Woronowicz (\cite{wor2}, Proposition 5.4). He proved that every irreducible representation enters the right regular representation with the multiplicity equal to its dimension. The result is certainly known for classical compact groups (see for instance \cite{ross}, Section 27.49). A. Maes and A. Van Daele (\cite{mvan}, Theorem 6.7) gave a simpler proof to the fact that every irreducible representation is contained in the right regular representation.  The aim of this paper is to give a proof of Woronowicz's result, including the multiplicities, which is an adaptation of the corresponding result for classical compact groups and is different from both Woronowicz's and Maes and Van Daele's proofs. The proof we give is based on the structure of $C^{\ast}(G)$ and its commutant and provides a concrete decomposition of the right regular representation in irreducible components.
\begin{proposition}[\cite{wor2}, Proposition 5.4]\label{thm} The right regular representation $v$ contains all irreducible representations of $G$. Each irreducible representation enters $v$ with the multiplicity equal to its dimension, i.e.
$v=\underset{\alpha \in \widehat{G}}{\sum}\oplus d_{\alpha}u^{\alpha}$.
\end{proposition}

\begin{proof} 

\ Let $\alpha \in \widehat{G}$. We will prove that $H_{h}^{k}(\overline{u^{\alpha}})$, $1\leqslant k\leqslant d_{\alpha}$ are invariant for $V$ and the restriction of $V$ to any of these subspaces is equivalent to $u^{\overline{\alpha}}$. Since $\forall \beta \in \widehat G$, $\exists \alpha \in \widehat G$ such that $u^{\overline{\alpha}}\sim u^{\beta}$ then the conclusion will follow.

\ Using Lemma \ref{V} with $a=\overset{d_{\alpha}}{\underset{i=1}{\sum}}\lambda_{i}{u^{\alpha}_{ki}}^{\ast}$ we obtain:
\begin{align}
V(\pi_{h}(\overset{d_{\alpha}}{\underset{i=1}{\sum}}\lambda_{i}{u^{\alpha}_{ki}}^{\ast})\xi_{h}\otimes \eta)&=(\pi_{h}\otimes \pi_{h})(\Delta(\overset{d_{\alpha}}{\underset{i=1}{\sum}}\lambda_{i}{u^{\alpha}_{ki}}^{\ast}))(\xi_{h}\otimes \eta)\notag \\
&=\overset{d_{\alpha}}{\underset{i,j=1}{\sum}}\lambda_{i}\pi_{h}({u^{\alpha}_{kj}}^{\ast})\xi_{h}\otimes \pi_{h}({u^{\alpha}_{ji}}^{\ast})\xi_{h}\in H_{h}^{k}(\overline{u^{\alpha}})\otimes H_{h}.\notag
\end{align}
\ Therefore, with $p_{\alpha}^{k}$ the orthogonal projection onto $H_{h}^{k}(\overline{u^{\alpha}})$, we obtain
\[V(p_{\alpha}^{k}\otimes 1)=(p_{\alpha}^{k}\otimes 1)V(p_{\alpha}^{k}\otimes 1)\] and hence (since $p_{\alpha}^{k}$ are orthogonal projections by Theorem \ref{comp})
\[ V(p_{\alpha}^{k}\otimes 1)=(p_{\alpha}^{k}\otimes 1)V.\]

\ Hence $p_{\alpha}^{k}$ is intertwining for $V$, so by Theorem \ref{inter}, $p_{\alpha}^{k}\in (C^{\ast}(G))'$. Since $p_{\alpha}\in C^{\ast}(G)$ it follows that $p_{\alpha}^{k}\in (C^{\ast}(G)p_{\alpha})'$. 

\ We prove next that the restriction of $v$ to any $H_{h}^{k}(\overline{u^{\alpha}})$ is equivalent to $\overline{u^{\alpha}}$ and therefore equivalent to $u^{\overline{\alpha}}$. This restriction is given by $v(p_{\alpha}^{k}\otimes 1)$.

\ Let $k_{0}\in \{1,...,d_{\alpha}\}$. We show that $V(p_{\alpha}^{k_{0}}\otimes 1)\sim \overline{U^{\alpha}}$, where $\overline{U^{\alpha}}=(id\otimes \pi_{h})\overline{u^{\alpha}}$, by proving that $\exists S\in B(\overline{H_{\alpha}},H_{h})$ such that
\begin{equation}\label{eqn}
V(p_{\alpha}^{k_{0}}\otimes 1)(S\otimes 1)=(S\otimes 1)\overline{U^{\alpha}}.
\end{equation}
\ Define $S:\overline{H_{\alpha}}\to H_{h}$ by:
\[ S\xi_{i}=\pi_{h}({u^{\alpha}_{k_{0}i}}^{\ast})\xi_{h}\in H_{h}^{k}(\overline{u^{\alpha}})\subset H_{h}\]
where $\{\xi_{i}\mid \:1\leqslant i\leqslant d_{\alpha}\}$ is an orthonormal basis in $\overline{H_{\alpha}}$, the Hilbert space of the representation $\overline{u^{\alpha}}$. Note that $\overline{H_{\alpha}}$ coincides with $H_{\overline{\alpha}}$. 

\ Then, with $\xi_{i_{0}}$ in the basis above and $\eta \in H_{h}$ we have
\begin{align}
V(p_{\alpha}^{k_{0}}\otimes 1)(S\otimes 1)(\xi_{i_{0}}\otimes \eta)&=V(p_{\alpha}^{k_{0}}\otimes 1)(\pi_{h}({u^{\alpha}_{k_{0}i_{0}}}^{\ast})\xi_{h}\otimes \eta)\notag\\
&=V(\pi_{h}({u^{\alpha}_{k_{0}i_{0}}}^{\ast})\xi_{h}\otimes \eta)\notag\\
&=(\pi_{h}\otimes \pi_{h})(\Delta({u^{\alpha}_{k_{0}i_{0}}}^{\ast}))(\xi_{h}\otimes \eta)\notag \\
&=\overset{d_{\alpha}}{\underset{j=1}{\sum}}\pi_{h}({u^{\alpha}_{k_{0}j}}^{\ast})\xi_{h}\otimes\pi_{h}({u^{\alpha}_{ji_{0}}}^{\ast})\eta, \notag
\end{align}
\ and
\begin{align}
(S\otimes 1)\overline{U^{\alpha}}(\xi_{i_{0}}\otimes \eta)&=(S\otimes 1)(id\otimes \pi_{h})\overline{u^{\alpha}}(\xi_{i_{0}}\otimes \eta)\notag \\
&=(S\otimes 1)(id\otimes \pi_{h})(\overset{d_{\alpha}}{\underset{i,j=1}{\sum}} m_{ij}\otimes {u^{\alpha}_{ij}}^{\ast})(\xi_{i_{0}}\otimes \eta)\notag \\
&=(S\otimes 1)(\overset{d_{\alpha}}{\underset{i,j=1}{\sum}} m_{ij}\xi_{i_{0}}\otimes \pi_{h}({u^{\alpha}_{ij}}^{\ast})\eta)\notag\\
&=(S\otimes 1)(\overset{d_{\alpha}}{\underset{i,j=1}{\sum}} \delta_{ji_{0}}\xi_{i}\otimes \pi_{h}({u^{\alpha}_{ij}}^{\ast})\eta)\notag \\
&=\overset{d_{\alpha}}{\underset{i=1}{\sum}}\pi_{h}({u^{\alpha}_{k_{0}i}}^{\ast})\xi_{h}\otimes\pi_{h}({u^{\alpha}_{ii_{0}}}^{\ast})\eta.\notag
\end{align}
\ Hence $V(p_{\alpha}^{k_{0}}\otimes 1)(S\otimes 1)=(S\otimes 1)\overline{U^{\alpha}}$. From this relation the conclusion will follow using the following arguments: since $v(p_{\alpha}\otimes 1)$ is a finite dimensional representation then $v(p_{\alpha}\otimes 1)\in B(\overline{H_{\alpha}})\otimes \mathcal{A}$. Since $h$ is faithful on $\mathcal{A}$, \eqref{eqn} above implies $v(p_{\alpha}^{k}\otimes 1)\sim \overline{u^{\alpha}}$ for all $k$ and, since $\overline{u^{\alpha}}\sim u^{\overline{\alpha}}$, the conclusion follows.
\end{proof}

\

\end{document}